\documentclass[10pt]{amsart}
\usepackage{mathptmx}
\usepackage{amsmath}
\usepackage{amssymb}
\usepackage{array}
\usepackage{geometry}
\usepackage[bookmarks=true,colorlinks=true, pdfstartview=FitV, linkcolor=black, citecolor=blue, urlcolor=black]{hyperref}

\usepackage{color}
\definecolor{DarkRed}{rgb}{0.55,.00,0.2}
\definecolor{DarkGrey}{rgb}{0.35,.35,0.35}

\theoremstyle{definition}

\theoremstyle{remark}

\numberwithin{equation}{section}

\begin{document}

\title{Lebedev's type   index transforms \\ \vspace{0,10cm}   with the modified Bessel functions}

\author{S. Yakubovich}
\address{Department of Mathematics, Faculty of Sciences,  University of Porto,  Campo Alegre str.,  687; 4169-007 Porto,  Portugal}
\email{ syakubov@fc.up.pt}



\date{\today}
\maketitle

\markboth{\rm \centerline{ S.  Yakubovich}}{}
\markright{\rm \centerline{Lebedev's Type  Index Transforms}}

\begin{abstract}  New index transforms of the Lebedev type are investigated. It involves the real part of the product of the modified Bessel functions as the kernel.   The boundedness and  invertibility are examined  for these operators    in the Lebesgue weighted spaces.  Inversion theorems are proved.   Important  particular cases are exhibited.   The results are applied to   solve an initial   value problem for the fourth  order PDE, involving the Laplacian.   Finally,  it is shown that  the same PDE has another  fundamental solution, which is  associated with the generalized Lebedev index transform, involving the square of the modulus of Macdonald's function, recently considered  by the author. 
\end{abstract}

\section{Introduction and preliminary results}

Let $\alpha \in \mathbb{R}$.  The objects  of this paper are the following  index transforms \cite{yak},\ \cite{yal} 
$$(F_\alpha f) (\tau) = {2\sqrt\pi\over \cosh(\pi\tau)}  \int_0^\infty {\rm Re} \left[ K_{\alpha + i\tau}(\sqrt x)\  I_{\alpha - i\tau}(\sqrt x)\right] f(x)dx,   \quad    \tau \in \mathbb{R}, \eqno(1.1)$$
$$( G_\alpha g) (x) =  2\sqrt\pi \int_{-\infty}^\infty {\rm Re} \left[ K_{\alpha + i\tau}(\sqrt x)\  I_{\alpha - i\tau}(\sqrt x)\right] {g(\tau)\over \cosh(\pi\tau)} d\tau,   \quad   x  \in \mathbb{R}_+, \eqno(1.2)$$
where $i$ is the imaginary unit,  ${\rm Re} $ denotes the real part of the complex-valued function and $K_\mu(z), I_\mu(z)$  \cite{erd}, Vol. II  are modified Bessel functions,  satisfying  the differential equation
$$  z^2{d^2u\over dz^2}  + z{du\over dz} - (z^2+\mu^2)u = 0.\eqno(1.3)$$
The asymptotic behaviour at infinity is given by the formulas 
$$ K_\mu(z) = \left( \frac{\pi}{2z} \right)^{1/2} e^{-z} [1+
O(1/z)], \qquad z \to \infty,\eqno(1.4)$$
$$ I_\mu(z) = \left( \frac{1} {2 \pi z} \right)^{1/2} e^{z} [1+
O(1/z)], \qquad z \to \infty.\eqno(1.5)$$
Near zero we have, correspondingly,  the relations 
$$K_\mu(z) = O( z^{- |{\rm Re} \mu|} ),\ \mu \neq 0,\quad K_0(z)= O( \log z ),\ z \to 0,\eqno(1.6)$$
$$I_\mu(z) = O( z^{ |{\rm Re} \mu|} ),\ z \to 0.\eqno(1.7)$$ 

The product of the modified Bessel functions can be represented by the following integrals (see relations (2.12.14.1) in \cite{prud}, Vol. II and  (8.4.23.25) in \cite{prud}, Vol. III)
 $$ {\rm Re} \left[ K_{\alpha + i\tau}(\sqrt x)\  I_{\alpha - i\tau}(\sqrt x)\right] = \int_0^\infty J_{2\alpha} \left( 2\sqrt x \ \sinh y\right) \cos(2y \tau ) dy,\ x > 0,\eqno(1.8)$$

$$ {\rm Re} \left[ K_{\alpha + i\tau}(\sqrt x)\  I_{\alpha - i\tau}(\sqrt x)\right] = {\cosh (\pi\tau) \over 4 \pi \sqrt\pi i } 
\int_{\gamma-i\infty}^{\gamma+i\infty} \Gamma(s+i\tau) \Gamma (s-i\tau) \frac{\Gamma (s+\alpha) \Gamma (1/2 -s)}
{ \Gamma(1+\alpha -s)  \Gamma(s)} x^{-s} ds,\  x >0,\eqno(1.9)$$
where $\max ( - \alpha, 0 ) < \gamma <  1/2$ and $J_\mu(z)$ is the Bessel function of the first kind \cite{erd}, Vol. II.   The latter integral is a key ingredient to derive the differential equation for the kernel of the index transforms (1.1), (1.2). In fact, denoting by 
$$\Phi_{\alpha, \tau}  (x)=  {2\sqrt\pi\over \cosh(\pi\tau)}   {\rm Re} \left[ K_{\alpha + i\tau}(\sqrt x)\  I_{\alpha - i\tau}(\sqrt x)\right],\eqno(1.10)$$
we have 

{\bf Lemma 1}. {\it The kernel  $\Phi_{\alpha, \tau}  (x) $ is a fundamental solution of  the following  fourth order differential equation with variable coefficients}
 $$ x^3 {d^4 \Phi_{\alpha, \tau} \over dx^4} +  6x^2 \  {d^3 \Phi_{\alpha, \tau} \over dx^3} + x \left( 7+\tau^2-\alpha^2 -x \right)  {d^2 \Phi_{\alpha, \tau} \over dx^2}  +  \left( 1+\tau^2-\alpha^2 - {5\over 2} \ x \right)  {d \Phi_{\alpha, \tau} \over dx}$$
 $$- \left( {(\alpha \tau)^2\over x } + {1\over 2} \right)\  \Phi_{\alpha, \tau} = 0.  \eqno(1.11)$$

\begin{proof}  Indeed, recalling the representation (1.9), we appeal to the Stirling asymptotic formula for the gamma-function \cite{erd}, Vol. I to write for each $\tau \in \mathbb{R}$ and $s= \gamma +it$
$$\Gamma(s+i\tau) \Gamma (s-i\tau) \frac{\Gamma (s+\alpha) \Gamma (1/2 -s)}
{ \Gamma(1+\alpha -s)  \Gamma(s)} = O\left( e^{-\pi |t|} |t|^{ 2\gamma - 3/2}\right), \quad |t| \to \infty.\eqno(1.12)  $$
This means that the repeated differentiation with respect to $x$ under the integral sign in (1.9) is allowed,  and with the use of the reduction formula for the gamma-function \cite{erd}, Vol. I we obtain 
$$\left( x{d\over dx} \right)^2 \Phi_{\alpha, \tau}  (x) =  {1\over 2\pi i} \int_{\gamma-i\infty}^{\gamma+i\infty} \Gamma(s+i\tau) \Gamma (s-i\tau) \frac{ s^2  \  \Gamma (s+\alpha) \Gamma (1/2 -s)}  { \Gamma(1+\alpha -s)  \Gamma(s)} x^{-s} ds$$
$$= - \tau^2  \Phi_{\alpha, \tau}  (x) + {1\over 2\pi i} \int_{\gamma-i\infty}^{\gamma+i\infty} \Gamma(1+ s+i\tau) \Gamma (1+ s-i\tau) \frac{ \Gamma (s+\alpha) \Gamma (1/2 -s)}  { \Gamma(1+\alpha -s)  \Gamma(s)} x^{-s} ds$$

$$=  - \tau^2  \Phi_{\alpha, \tau}  (x) + {1\over 2\pi i} \int_{1+ \gamma-i\infty}^{1+ \gamma+i\infty} \Gamma(s+i\tau) \Gamma ( s-i\tau) \frac{ \Gamma (s-1 +\alpha) \Gamma (3/2 -s)}  { \Gamma(2+\alpha -s)  \Gamma(s-1)} x^{1 -s} ds$$
$$=  - \tau^2  \Phi_{\alpha, \tau}  (x) + {1\over 2\pi i} \int_{1+ \gamma-i\infty}^{1+ \gamma+i\infty} \Gamma(s+i\tau) \Gamma ( s-i\tau) \frac{ \Gamma (s +\alpha) \Gamma (1/2 -s) (1/2-s) (s-1) }  { (s-1 +\alpha) (1+\alpha -s)  \Gamma(1+\alpha -s)  \Gamma(s)} x^{1 -s} ds.$$
Hence, moving the contour to the left by Cauchy's theorem, we multiply the latter equality by $x^\alpha$ and differentiate with respect to $x$ again. Then  
$${d\over dx} \left[ x^\alpha \left( x{d\over dx} \right)^2 \Phi_{\alpha, \tau}  (x) \right] =  - \tau^2 {d\over dx} \left[ x^\alpha  \Phi_{\alpha, \tau}  (x) \right] $$
$$+ {1\over 2\pi i} \int_{ \gamma-i\infty}^{\gamma+i\infty} \Gamma(s+i\tau) \Gamma ( s-i\tau) \frac{ \Gamma (s +\alpha) \Gamma (1/2 -s) (1/2-s) (s-1) }  { (s-1 +\alpha)  \Gamma(1+\alpha -s)  \Gamma(s)} x^{\alpha -s} ds.$$
In a similar manner we continue to reduce the denominator of the integrand, multiplying by $x^{1-2\alpha}$ and fulfilling again the differentiation.  Hence,
$$ {d\over dx} \left[ x^{1-2\alpha}  {d\over dx} \left[ x^\alpha \left( x{d\over dx} \right)^2 \Phi_{\alpha, \tau}  (x) \right] \right] = 
 - \tau^2 {d\over dx} \left[ x^{1-2\alpha}  {d\over dx} \left[ x^\alpha  \Phi_{\alpha, \tau}  (x) \right] \right] $$
$$ -  {1\over 2\pi i} \int_{ \gamma-i\infty}^{\gamma+i\infty} \Gamma(s+i\tau) \Gamma ( s-i\tau) \frac{ \Gamma (s +\alpha) \Gamma (1/2 -s) (1/2-s) (s-1) }  { \Gamma(1+\alpha -s)  \Gamma(s)} x^{- \alpha -s} ds.$$
Now multiplying by $x^{\alpha}$ and accounting (1.9), we find 
$$ x^{\alpha} {d\over dx} \left[ x^{1-2\alpha}  {d\over dx} \left[ x^\alpha \left( x{d\over dx} \right)^2 \Phi_{\alpha, \tau}  (x) \right] \right] =   - \tau^2\  x^{\alpha}  {d\over dx} \left[ x^{1-2\alpha}  {d\over dx} \left[ x^\alpha  \Phi_{\alpha, \tau}  (x) \right] \right] $$
$$ +  {1\over 2\pi i}  {d\over dx} \int_{ \gamma-i\infty}^{\gamma+i\infty} \Gamma(s+i\tau) \Gamma ( s-i\tau) \frac{ \Gamma (s +\alpha) \Gamma (1/2 -s) (1/2-s) }  { \Gamma(1+\alpha -s)  \Gamma(s)} x^{1 -s} ds$$
$$=  - \tau^2\  x^{\alpha}  {d\over dx} \left[ x^{1-2\alpha}  {d\over dx} \left[ x^\alpha  \Phi_{\alpha, \tau}  (x) \right] \right] - {3\over 2} \  {d \over dx} \left[ x\  \Phi_{\alpha, \tau}  (x)\right]  +  {d^2 \over dx^2} \left[ x^2\  \Phi_{\alpha, \tau}  (x)\right].$$ 
Finally, fulfilling the differentiation, we end up with the equation (1.11), completing the proof of Lemma 1. 
\end{proof}

\section {Boundedness and invertibility properties for  the index transform (1.1)}

Our approach to examine  the boundedness and invertibility  properties of the introduced index transforms  is based on  the Mellin transform technique developed in \cite{yal} and extensive use of the Marichev  table for the Mellin transform in \cite{mar}, \cite{prud}, Vol. III.   Appealing to the classical Titchmarsh monograph \cite{tit},  the Mellin transform is defined, for instance, in  $L_{\nu, p}(\mathbb{R}_+),\ 1 < p \le 2$ by the integral  
$$f^*(s)= \int_0^\infty f(x) x^{s-1} dx,\eqno(2.1)$$
 where the convergence is understood  in mean with respect to the norm in $L_q(\nu- i\infty, \nu + i\infty),\   q=p/(p-1)$.   Moreover, the  Parseval equality holds for $f \in L_{\nu, p}(\mathbb{R}_+),\  g \in L_{1-\nu, q}(\mathbb{R}_+)$
$$\int_0^\infty f(x) g(x) dx= {1\over 2\pi i} \int_{\nu- i\infty}^{\nu+i\infty} f^*(s) g^*(1-s) ds.\eqno(2.2)$$
The inverse Mellin transform is given accordingly
 $$f(x)= {1\over 2\pi i}  \int_{\nu- i\infty}^{\nu+i\infty} f^*(s)  x^{-s} ds,\eqno(2.3)$$
where the integral converges in mean with respect to the norm  in   $L_{\nu, p}(\mathbb{R}_+)$
$$||f||_{\nu,p} = \left( \int_0^\infty  |f(x)|^p x^{\nu p-1} dx\right)^{1/p}.\eqno(2.4)$$
In particular, letting $\nu= 1/p$ we get the usual space $L_1(\mathbb{R}_+)$.  Further, denoting by $C (\mathbb{R})$ the space of bounded continuous functions, we prove the following result. 

{\bf Theorem 1.}   {\it  Let $\alpha > - 1/4$. The index transform  $(1.1)$  is well-defined as a  bounded operator $F_\alpha : L_{3/4, 1} \left(\mathbb{R}_+ \right) \to C (\mathbb{R})$.   Moreover,  if in addition $f \in L_{1-\nu, p}(\mathbb{R}_+),\  1 < p \le 2,\  \max ( - \alpha, 0 ) < \nu <  1/2$,  then 
$$(F_\alpha f)(\tau)=  {2 \sqrt\pi\over  \cosh(\pi\tau)} \int_0^\infty    K_{ i\tau}(\sqrt x) {\rm Re} \left[  I_{i\tau}(\sqrt x)\right] \varphi_\alpha  (x) dx,\eqno(2.5)$$
where the integral converges absolutely, 
$$\varphi_\alpha (x)= {1\over 2\pi i}  \int_{1-\nu- i\infty}^{1-\nu+i\infty}      \frac{\Gamma (1- s+\alpha) \Gamma (s)}
{ \Gamma(s+\alpha )  \Gamma(1- s)} \   f^*(s)   x^{-s} ds,\eqno(2.6)$$
and  integral $(2.6)$ converges in mean with respect to the norm in $L_{1-\nu, p}(\mathbb{R}_+)$ .}

\begin{proof}  Recalling the integral representation (1.8) of the index kernel in (1.1) and elementary inequality for the Bessel function of the first kind  $\sqrt x \left|J_\mu (x) \right| < C,\  {\rm Re\  \mu }  > -1/2,\ C >0$ is an absolute constant, 
we have the estimate (we will keep the same notation for  different  positive constants )
$$ \left| (F_\alpha  f) (\tau)\right| \le  {2\sqrt\pi\over \cosh(\pi\tau)}  \int_0^\infty \int_0^\infty \left| J_{2\alpha} \left( 2\sqrt x \ \sinh y\right) \right| |f(x) |  dy\ dx $$
$$\le C   \int_0^\infty {dy\over \sqrt {\sinh y} }  \int_0^\infty x^{-1/4}  |f(x) | \  dx  = C\   ||f||_{3/4, 1},\  \alpha > - {1\over 4}.$$
Hence via the absolute and uniform convergence it follows the continuity of $(F_\alpha f)(\tau)$ and the boundedness of the operator (1.1), namely
$$\sup_{\tau \in \mathbb{R}}   \left| (F_\alpha  f) (\tau)\right| \equiv || F_\alpha f||_ {C (\mathbb{R})} \le C ||f||_{3/4, 1}.$$
Further,  from the condition  $f \in L_{1-\nu, p}(\mathbb{R}_+)$, asymptotic behaviour (1.12) and the Parseval equality (2.2) we derive the representation
$$  (F_\alpha  f) (\tau) =  {1 \over 2 \pi i }  \int_{\nu -i\infty}^{\nu +i\infty} \Gamma(s+i\tau) \Gamma (s-i\tau) \frac{\Gamma (s+\alpha) \Gamma (1/2 -s)} { \Gamma(1+\alpha -s)  \Gamma(s)}  f^*(1-s) ds.\eqno(2.7)$$
Meanwhile, by virtue of the Stirling formula for the gamma-function 
$$ \frac{\Gamma (1- s+\alpha) \Gamma (s)}  { \Gamma(s+\alpha )  \Gamma(1- s)} = O(1),\  s= 1-\nu +it,\ |t| \to \infty.\eqno(2.8)$$
Therefore, employing relation (8.4.23.23) in \cite{prud}, Vol. III, we apply again the Parseval equality (2.2) to the right-hand side of (2.7).  This leads to the formula  (2.5), which is, in turn, the Lebedev index transform with the modified Bessel functions as the kernel \cite{square} of the function $\varphi_\alpha$  defined by (2.6).  Theorem 1 is proved. 

\end{proof}

The inversion formula for the Lebedev type transform (1.1) is established by

{\bf Theorem 2.}  {\it Let $\alpha > 1/4,\  f \in  L_{1-\nu, p}(\mathbb{R}_+)$, its  Mellin transform $f^*(s) \in L_1( 1-\nu-i\infty, 1-\nu +i\infty)$ and be analytic in the open strip ${\rm Re} s = 1- \nu  \in (0,  \alpha + 1)$.    Then under  the integrability condition $(F_\alpha f)(\tau) \in L_1(\mathbb{R}_+; \tau e^{2 \pi\tau} d\tau)$ the following inversion formula  holds for all $x >0$ 

$$ f(x)= -  {4 \over \pi \sqrt\pi }   \int_0^\infty  \left[  \alpha  \left[ {1\over 2 x}   +  {d\over dx} \left( {\rm Re} \  \left[ K _{\alpha + i\tau}(\sqrt x )\  I_{\alpha  - i\tau}(\sqrt x ) \right] \ \right)  - \sqrt x \   {\rm Re} \  \left[ K _{\alpha + i\tau}(\sqrt x )\  I^\prime_{\alpha  - i\tau}(\sqrt x ) \  \right]\right] \right.$$
$$\left. +     \tau   \  {d\over dx} \left( {\rm Im} \  \left[ K _{\alpha + i\tau}(\sqrt x )\  I_{\alpha  - i\tau}(\sqrt x ) \right] \  \right ) \right]\    \cosh(\pi\tau)  (F_\alpha f) (\tau) d\tau dy,\eqno(2.9)$$
where $\prime$ is the symbol for  derivative,  ${\rm Im} $ denotes the imaginary part of a complex -valued function and the corresponding  integral converges absolutely.}

\begin{proof}   Since $f \in  L_{1-\nu, p}(\mathbb{R}_+)$ we have from Theorem 86 in \cite{tit} that  $f^*(s) \in L_q(1-\nu-i\infty, 1-\nu +i\infty),\ q= p/(p-1)$.   Recalling (2.5), our goal now is to verify conditions  of the Lebedev expansion theorem in \cite{square} in order to prove the following inversion formula
$$\int_x^\infty  \varphi_\alpha (y) dy = {2\over \pi^2 \sqrt\pi} \int_0^\infty \tau \sinh(2\pi\tau) K_{i\tau}^2 (\sqrt x) (F_\alpha f) (\tau) d\tau,\ x >0.\eqno(2.10)$$
To do this,  we show that $\varphi_\alpha \in L_{3/4, 1} ((0,1)) \cap  L_{5/4, 1} ((1, \infty))$.   Indeed,  the integrability condition, estimate (2.8)  and analyticity of $f^*(s)$ in the strip ${\rm Re} s = 1- \nu  \in (0,  \alpha + 1)$ allow to to move the contour in (2.6) by Cauchy's theorem , keeping the same value.   Therefore, choosing $\nu \in (1/4, 1)$, we use the H\"{o}lder inequality to find
$$||\varphi_\alpha ||_{3/4,1} = \int_0^1 |\varphi_\alpha(x) |  x^{- 1/4} dx \le ||\varphi_\alpha ||_{1-\nu,p} \left( \int_0^1 x^{(\nu- 1/4)q-1} dx \right)^{1/q} = { ||\varphi_\alpha ||_{1-\nu,p} \over [(\nu-1/4) q] ^{1/q}} <\infty.$$
In the meantime, taking $\nu \in (- \alpha, -1/4)$, we have 
$$||\varphi_\alpha ||_{5/4,1} = \int_1^\infty  |\varphi_\alpha(x) |  x^{1/4} dx \le ||\varphi_\alpha ||_{1-\nu,p} \left( \int_1^\infty  x^{(\nu+ 1/4)q-1} dx \right)^{1/q} = { ||\varphi_\alpha ||_{1-\nu,p} \over [- (\nu+ 1/4) q] ^{1/q}} <\infty.$$
Therefore,  substituting (2.6) in the left-hand side of (2.10) and making the integration with respect to $y$ via Fubini's theorem and a simple substitution, we obtain 
$$-  {1\over 2\pi i}  \int_{-\nu- i\infty}^{-\nu+i\infty}      \frac{\Gamma (\alpha -s) \Gamma (s+1)}
{ \Gamma(1+ s+\alpha )  \Gamma(1- s)} \   f^*(1+s)   x^{-s} ds = {2\over \pi^2 \sqrt\pi} \int_0^\infty \tau \sinh(2\pi\tau) K_{i\tau}^2 (\sqrt x) (F_\alpha f) (\tau) d\tau.\eqno(2.11)$$
Further, taking the Mellin transform from both sides of (2.11), basing on the condition $(F_\alpha f)(\tau) \in \  L_1(\mathbb{R}_+; \\  \tau e^{2\pi\tau} d\tau)$ and the uniform estimate $|K_{i\tau} (\sqrt x)| \le K_0(\sqrt x) $ for the modified Bessel function, we change the order of integration due to the absolute and uniform convergence and appeal to relation (8.4.23.27) in \cite{prud}, Vol. III to calculate the inner integral with respect to $x$ in the right-hand side of the obtained equality. Hence,  recalling the reduction formula for the gamma-function, we end up with the equality
$$ s f^*(1+ s) = - {1\over \pi^2}   \frac { \Gamma(1+ s+\alpha )  \Gamma(1- s)} {\Gamma (\alpha -s) \Gamma (s+1/2)}
 \int_0^\infty \tau \sinh(2\pi\tau) \Gamma(s+i\tau) \Gamma (s-i\tau) (F_\alpha f) (\tau) d\tau.\eqno(2.12)$$
Then,  reciprocally,  via formula (2.3) and properties of the Mellin transform we find 
$$x f(x)=  {1\over \pi^2}   \int_0^\infty \tau \sinh(2\pi\tau) S(x,\tau) (F_\alpha f) (\tau) d\tau ,\eqno(2.13)$$
where
$$S(x,\tau) =   -  {1\over 2\pi i}  \int_x^\infty  \int_{\gamma - i\infty}^{\gamma +i\infty}    \Gamma(s+i\tau) \Gamma (s-i\tau) 
\frac { \Gamma(1+ s+\alpha )  \Gamma(1- s)} {\Gamma (\alpha -s) \Gamma (s+1/2)}  y^{-s-1} ds dy \eqno(2.14)$$
and $\gamma \in (0,1)$.  The kernel $S(y,\tau)$ can be calculated, recalling relation (8.4.23.25) in \cite{prud}, Vol. III and using  repeated differentiation under the integral sign.  In fact, via relation (8.4.23.25) we have that 
$$ {\sqrt\pi\over i \sinh (\pi\tau)} \left[ K_{\alpha + i\tau}(\sqrt x)\  I_{\alpha - i\tau}(\sqrt x) -  K_{\alpha - i\tau}(\sqrt x)\  I_{\alpha +  i\tau}(\sqrt x) \right] $$
$$= {1\over 2\pi i} \int_{\gamma - i\infty}^{\gamma +i\infty}    \Gamma(s+i\tau) \Gamma (s-i\tau) 
\frac { \Gamma( s+\alpha )  \Gamma(1- s)} {\Gamma (1+ \alpha -s) \Gamma (s+1/2)}  x^{-s} ds,\  x >0,\  \tau \in \mathbb{R} \backslash \{0\}.\eqno(2.15) $$
Hence, involving the repeated differentiation and the reduction formula for the gamma-function, it is not difficult to verify the equalities 
$$ -  {1\over 2\pi i}   \int_{\gamma - i\infty}^{\gamma +i\infty}    \Gamma(s+i\tau) \Gamma (s-i\tau) 
\frac { \Gamma(1+ s+\alpha )  \Gamma(1- s)} {\Gamma (\alpha -s) \Gamma (s+1/2)}  x^{-s-1} ds$$
$$=   {\sqrt\pi\over i \sinh (\pi\tau)}  x^\alpha {d\over dx } \ x^{1-2\alpha}  {d\over dx } \ x^{\alpha} 
\left[ K_{\alpha + i\tau}(\sqrt x)\  I_{\alpha - i\tau}(\sqrt x) -  K_{\alpha - i\tau}(\sqrt x)\  I_{\alpha +  i\tau}(\sqrt x) \right] $$
$$= {\sqrt\pi\over i \sinh (\pi\tau)} \left[   {d\over dx } \ x  {d\over dx }  \left[ K_{\alpha + i\tau}(\sqrt x)\  I_{\alpha - i\tau}(\sqrt x) -  K_{\alpha - i\tau}(\sqrt x)\  I_{\alpha +  i\tau}(\sqrt x) \right]  \right. $$
$$\left.  - {\alpha^2\over x} \left[ K_{\alpha + i\tau}(\sqrt x)\  I_{\alpha - i\tau}(\sqrt x) -  K_{\alpha - i\tau}(\sqrt x)\  I_{\alpha +  i\tau}(\sqrt x) \right] \right].$$
Therefore, talking into account the asymptotic behaviour of the function
$$x  {d\over dx }  \left[ K_{\alpha + i\tau}(\sqrt x)\  I_{\alpha - i\tau}(\sqrt x) -  K_{\alpha - i\tau}(\sqrt x)\  I_{\alpha +  i\tau}(\sqrt x) \right] = o(1),\ x \to +\infty,$$
which can be established, for instance, from the integral representation (2.15) and Stirling asymptotic formula for the gamma-function, we return to (2.14)  to obtain 
$$S(x,\tau) = -  { \sqrt\pi  \over i \sinh (\pi\tau)} \left[ 2\alpha^2   \int_{\sqrt x}^\infty   \left[ K_{\alpha + i\tau}(y )\  I_{\alpha - i\tau}(y ) -  K_{\alpha - i\tau}(y)\  I_{\alpha +  i\tau}(y ) \right] {dy\over y} \right.$$ 
$$\left.   +  x   {d\over dx }  \left[ K_{\alpha + i\tau}(\sqrt x)\  I_{\alpha - i\tau}(\sqrt x) -  K_{\alpha - i\tau}(\sqrt x)\  I_{\alpha +  i\tau}(\sqrt x) \right] \right].\eqno(2.16)$$
However, the integral in (2.16) can be treated, employing relations (1.12.4.4),  (2.16.28.3) in \cite{prud}, Vol. II and asymptotic formulae (1.4), (1.5), (1.6), (1.7) for the modified Bessel functions.  Then we derive
$$\int_{\sqrt x}^\infty   \left[ K_{\alpha + i\tau}(y )\  I_{\alpha - i\tau}(y ) -  K_{\alpha - i\tau}(y)\  I_{\alpha +  i\tau}(y ) \right] {dy\over y} $$
$$  = \lim_{\varepsilon \to 0+} \left(  \int_{0 }^\infty  -  \int_{0}^{\sqrt x} \right)    \left[ K_{\alpha + i\tau}(y )\  I_{\alpha +\varepsilon - i\tau}(y ) -  K_{\alpha - i\tau}(y)\  I_{\alpha +\varepsilon  +  i\tau}(y ) \right] {dy\over y}$$
$$= {i\over 2\alpha\tau} -   \lim_{\varepsilon \to 0+}  \int_{0}^{\sqrt x}  \left[ K_{\alpha + i\tau}(y )\  I_{\alpha +\varepsilon - i\tau}(y ) -  K_{\alpha - i\tau}(y)\  I_{\alpha +\varepsilon  +  i\tau}(y ) \right] {dy\over y}$$
$$=  {i\over 2\alpha\tau} +   {i \sqrt x\over 4\alpha \tau} \left[ K^\prime _{\alpha + i\tau}(\sqrt x )\  I_{\alpha  - i\tau}(\sqrt x )-
K _{\alpha + i\tau}(\sqrt x )\  I^\prime _{\alpha  - i\tau}(\sqrt x )\right] $$
$$+ {i \sqrt x\over 4\alpha \tau} \left[ K^\prime _{\alpha - i\tau}(\sqrt x )\  I_{\alpha  + i\tau}(\sqrt x )-
K _{\alpha - i\tau}(\sqrt x )\  I^\prime _{\alpha  + i\tau}(\sqrt x )\right] $$
$$=   {i\over 2\alpha\tau} +   {i  x\over 2\alpha \tau} {d\over dx} \left[ K _{\alpha + i\tau}(\sqrt x )\  I_{\alpha  - i\tau}(\sqrt x )+ K _{\alpha - i\tau}(\sqrt x )\  I_{\alpha  + i\tau}(\sqrt x ) \right]$$
$$ -  {i \sqrt x\over 2\alpha \tau} \left[ K _{\alpha + i\tau}(\sqrt x )\  I^\prime _{\alpha  - i\tau}(\sqrt x ) + K _{\alpha - i\tau}(\sqrt x )\  I^\prime _{\alpha  + i\tau}(\sqrt x )\right].$$
Hence, combining with (2.16),  we find
$$ S(x,\tau) =  -  {  \alpha \sqrt\pi \over  \tau \sinh (\pi\tau)} \left[ 1  + 2x {d\over dx} \left( {\rm Re} \  \left[ K _{\alpha + i\tau}(\sqrt x )\  I_{\alpha  - i\tau}(\sqrt x ) \right] \ \right)  - 2\sqrt x \   {\rm Re} \  \left[ K _{\alpha + i\tau}(\sqrt x )\  I^\prime_{\alpha  - i\tau}(\sqrt x ) \  \right]\right]$$
$$-   {  2x \ \sqrt\pi \over  \sinh (\pi\tau)}  {d\over dx} \left( {\rm Im} \  \left[ K _{\alpha + i\tau}(\sqrt x )\  I_{\alpha  - i\tau}(\sqrt x ) \right] \  \right).$$
Substituting this expression of $S(x,\tau) $ into (2.13),  we come up with inversion formula (2.9), completing the proof of Theorem 2. 

\end{proof}

{\bf Remark 1}.  Letting formally $\alpha =0$ in (2.9), we appeal to the relation (cf. \cite{erd}, Vol. II) for the Macdonald function 
$$K_{i\tau} (\sqrt x) = {\pi\over \sinh(\pi\tau) }  {\rm Im} \  \left[ I_{- i\tau}(\sqrt x ) \right],\eqno(2.17)$$
and making simple substitutions,  we arrive at the Lebedev inversion formula (2.10), where $\varphi_0(x) = f(x)$ via (2.3) .

\section{The index transform (1.2)} 

In this section we investigate the boundedness and invertibility properties for the Lebedev type transform  (1.2).

{\bf Theorem 3.}  {\it  Let $\alpha > - 1/4,\ g  \in L_1(\mathbb{R};\  [\cosh(\pi\tau)]^{-1}  d\tau)$. Then $x^{1/4} (G_\alpha g)(x)$ is bounded continuous on $\mathbb{R}_+$ and   
$$\sup_{x >0}  \   x^{1/4} \left| (G_\alpha g)(x) \right|   \le C  ||g ||_{L_1(\mathbb{R};\  [\cosh(\pi\tau)]^{-1}  d\tau)}.\eqno(3.1)$$
Besides,   if $g\in L_1(\mathbb{R})$ and $(G_\alpha g) (x) \in L_{\nu,1}((0,1)),\    \max ( - \alpha, 0 )   < \nu < 1/2$, then for all $y  >0$ }
$${1\over 2\pi i}   \int_{\nu -i\infty}^{\nu  +i\infty} \frac{\Gamma(s)\Gamma (1+\alpha -s)}{\Gamma(s+\alpha)}  (G_\alpha g)^*(s)  y^{ -s} ds\  =  \sqrt \pi  \   \int_{-\infty}^\infty  e^{y/2} \   K_{i\tau} \left({y\over 2} \right)  {g(\tau) \over \cosh(\pi\tau) } d\tau. \eqno(3.2)$$

\begin{proof}   Doing similarly as in the proof of Theorem 1, we employ (1.8) to have an immediate estimate
$$ \left| (G_\alpha g)(x) \right| \le C\  x^{-1/4} \int_{-\infty}^\infty {|g(\tau)| \over \cosh(\pi\tau)} d\tau= C  ||g ||_{L_1(\mathbb{R};\  [\cosh(\pi\tau)]^{-1}  d\tau)},$$
which yields (3.1).  In order to prove the equality (3.2),  we use the uniform estimate for the kernel (1.10), which can be obtained from the Mellin-Barnes representation (1.9). Indeed, by   definition  of the Euler beta-function \cite{erd}, Vol. I and Stirling asymptotic formula for the gamma-function we have for $x >0$
$$\left| \Phi_{\alpha, \tau}  (x) \right| \le {x^{-\gamma}  \over 2 \pi }  \int_{\gamma-i\infty}^{\gamma+i\infty} \left| \Gamma(s+i\tau) \Gamma (s-i\tau) \frac{\Gamma (s+\alpha) \Gamma (1/2 -s)}
{ \Gamma(1+\alpha -s)  \Gamma(s)} \right| ds $$
$$ \le  {x^{-\gamma}  B(\gamma,\gamma) \over 2 \pi }  \int_{\gamma-i\infty}^{\gamma+i\infty} \left| \Gamma(2 s) \frac{\Gamma (s+\alpha) \Gamma (1/2 -s)} { \Gamma(1+\alpha -s)  \Gamma(s)} \right| ds = C\  x^{-\gamma} \eqno(3.3)$$
where $\max ( - \alpha, 0 ) < \gamma <  1/2$.  Hence applying the Mellin transform (2.1) to both sides of (1.2),  we change the order of integration in the right-hand side of the obtained equality by Fubini's theorem.  The Mellin transform of its left-hand side exists under the condition $(G_\alpha g) (x) \in L_{\nu,1}((0,1)),\    \max ( - \alpha, 0 )   < \nu < 1/2$ and estimate (3.3), which guarantees the integrability of $(G_\alpha g) (x)$ over $(1,\infty)$. Therefore we end up with the equality 
$$\frac{\Gamma(s)\Gamma (1+\alpha -s)}{\Gamma(s+\alpha)}  (G_\alpha g)^*(s)  =  \Gamma(1/2-s) \int_{-\infty}^\infty    g(\tau) \Gamma(s+ i\tau)\Gamma(s-i\tau)  d\tau.$$
Finally, the inverse Mellin transform (2.3) and relation (8.4.23.5) in \cite{prud}, Vol. III will drive us to the equality (3.2), completing the proof of Theorem 3. 

\end{proof} 

Now we are ready to prove the inversion formula for the index transform (1.2).

{\bf Theorem 4}.  {\it  Let $g(z/i)$ be an even analytic function in the strip $D= \left\{ z \in \mathbb{C}: \ |{\rm Re} z | < \beta < 1/2\right\} ,\  g(0)=g^\prime (0)=0$ and $g(z/i)$ is absolutely  integrable over any vertical line in  $D$.   Then under conditions of Theorem 3  for all $x \in \mathbb{R}$ the  following inversion formula holds 
$$  g(x) = {\cosh (\pi x)\over \pi\sqrt\pi } \lim_{\varepsilon \to 0}  \int_0^\infty   \left[  y^\alpha \left| {\Gamma (\varepsilon-1-\alpha +ix) \over \Gamma (ix)}\right| ^2  \frac{ \Gamma(1+\alpha) }{ \Gamma(1+2\alpha)\Gamma(\varepsilon- 1/2-\alpha)} \right.$$
$$\times \  {}_2F_3 \left( 1+\alpha,\ {3\over 2}+\alpha-\varepsilon;\ 1+ 2\alpha,\ 2+\alpha -\varepsilon- ix,\  2+\alpha -\varepsilon+  ix ;\ y \right) $$
$$+ { y^{\varepsilon-1}  \over \sqrt\pi} {\rm Re } \left[ \left({y\over 4}\right)^{ix} \ \frac{  \Gamma(\varepsilon+ ix)} { \Gamma(ix)}  \  \Gamma(1- \varepsilon+\alpha -ix) \right. $$
$$\left. \left. \times  \  {}_2F_3 \left( \varepsilon+ix,\  {1\over 2}+ ix ;\ 1+ 2 ix,\  \varepsilon- \alpha + ix,\  \varepsilon +\alpha +  ix ;\ y \right) \right] \right] (G_\alpha g) (y) dy,\eqno(3.4)$$
where ${}_2F_3 (a_1,\ a_2; \ b_1,\ b_2,\ b_3; z)$ is the generalized hypergeometric function \cite{erd}, Vol. I}. 

\begin{proof}    In fact,  multiplying  both sides of (3.2) by $ e^{-y/2} K_{ix} \left({y/2} \right) y^{\varepsilon -1}$ for some positive $\varepsilon \in (0,1)$ we  integrate with respect to $y$ over $(0, \infty)$.  Hence,  since under conditions of the theorem $(G_\alpha g)^*(s)$ is bounded, we  change  the order of integration in the left-hand side of the obtained equality and   appeal  to the relation (8.4.23.3) in \cite{prud}, Vol. III.   Therefore,   for $\nu \in (0, \min( \varepsilon,\ 1+\alpha) )$
$$ {1\over 2\pi i}   \int_{\nu -i\infty}^{\nu  +i\infty} \frac {\Gamma(\varepsilon- s+ ix)\Gamma(\varepsilon- s-ix) \Gamma(s)\Gamma (1+\alpha -s)}{ \Gamma(1/2+ \varepsilon -s) \Gamma(s+\alpha)}  (G_\alpha g)^*(s) ds$$
$$ =  \int_0^\infty  K_{ix} \left({y\over 2} \right) y^{\varepsilon -1} \int_{-\infty}^\infty  K_{i\tau} \left({y\over 2} \right)  {g(\tau) \over \cosh(\pi\tau) } d\tau dy. \eqno(3.5)$$

In the meantime,  the right-hand side of (3.5) can be treated, taking into account   the evenness of $g$ and representation (2.17) for the Macdonald function.  Indeed, we have  
$$ \int_0^\infty  K_{ix} \left({y\over 2} \right) y^{\varepsilon -1} \int_{-\infty}^\infty  K_{i\tau} \left({y\over 2} \right)  {g(\tau) \over \cosh(\pi\tau) } d\tau dy$$

$$= 2 \pi i  \int_0^\infty  K_{ix} \left({y\over 2} \right) y^{\varepsilon -1} \int_{-i\infty}^{i\infty}   I_{ z} \left({y\over 2} \right)  {g(z/i) \over \sin (2\pi z) } dz\  dy. \eqno(3.6)$$
On the other hand, according to our assumption $g(z/i)$ is analytic in the vertical  strip $0\le  {\rm Re}  z < \beta < 1/2$,  $g(0)=g^\prime (0)=0$ and integrable  in the strip.  Hence,  appealing to the inequality for the modified Bessel   function of the first  kind  (see \cite{yal}, p. 93)
 $$|I_z(y)| \le I_{  {\rm Re} z} (y) \  e^{\pi |{\rm Im} z|/2},\   0< {\rm Re} z < \beta,$$
one can move the contour to the right in the latter integral in (3.6). Then 

$$2 \pi i  \ \int_0^\infty  K_{ix} \left({y\over 2} \right) y^{\varepsilon -1} \int_{-i\infty}^{i\infty}   I_{ z} \left({y\over 2} \right)  {g(z/i) \over \sin (2\pi z) } dz\  dy$$

$$= 2 \pi i    \int_0^\infty  K_{ix} \left({y\over 2} \right) y^{\varepsilon -1} \int_{\beta -i\infty}^{\beta + i\infty}   I_{ z} \left({y\over 2} \right)  {g(z/i) \over \sin (2\pi z) } dz\  dy.$$
Now ${\rm Re} z >0$,  and  it is possible to pass to the limit under the integral sign when $\varepsilon \to 0$ and to change the order of integration due to the absolute and uniform convergence.  Recalling the relation (2.16.28.3) in \cite{prud}, Vol. II, we find 
$$\lim_{\varepsilon \to 0}  2 \pi i   \int_0^\infty  K_{ix} \left({y\over 2} \right) y^{\varepsilon -1} \int_{-i\infty}^{i\infty}   I_{ z} \left({y\over 2} \right)  {g(z/i) \over \sin (2\pi z) } dz\  dy$$

$$=    2 \pi i  \ \int_{\beta -i\infty}^{\beta  + i\infty}   {g(z/i) \over (x^2+ z^2) \sin (2\pi z) } dz =  \pi i 
 \left( \int_{-\beta +i\infty}^{- \beta - i\infty}   +   \int_{\beta  -i\infty}^{ \beta +  i\infty}   \right)  {  g(z/i) \  dz \over (z-ix) \  z \sin(2\pi z)}. \eqno(3.7)$$
Hence, using the Cauchy formula in the right-hand side of the latter equality in (3.7) under conditions of the theorem, we derive 
$$\lim_{\varepsilon \to 0}  2 \pi i  \ \int_0^\infty  K_{ix} \left({y\over 2} \right) y^{\varepsilon -1} \int_{-i\infty}^{i\infty}   I_{ z} \left({y\over 2} \right)  {g(z/i) \over \sin (2\pi z) } dz\  dy =  { 2\pi^{2} \  g(x) \over  x\sinh (2\pi x)} ,\quad x \in \mathbb{R} \backslash \{0\}.\eqno(3.8)$$
Thus, returning to (3.5), employing the Parseval identity (2.2)  and passing to the limit when $\varepsilon \to 0$, we come up with the equality 
$$ { 2\pi^{2} \  g(x) \over  x\sinh (2\pi x)} = \lim_{\varepsilon \to 0}  \int_0^\infty S_\varepsilon (x,  y)  (G_\alpha g) (y) dy,\eqno(3.9)$$
where
$$S_\varepsilon (x,  y) =  {1\over 2\pi i}   \int_{1-\nu -i\infty}^{1- \nu  +i\infty} \frac {\Gamma(\varepsilon-1 + s+ ix)\Gamma(\varepsilon -1+ s-ix) \Gamma(s+\alpha )\Gamma (1 -s)}{ \Gamma(\varepsilon -1/2 + s) \Gamma(1- s+\alpha)}  y^{-s} ds.\eqno(3.10)$$
Meanwhile,  integral (3.10) can be calculated with the use of Slater's theorem \cite{mar} in terms of the generalized hypergeometric functions ${}_2F_3$.  Namely, it involves the left-hand simple poles of the gamma-functions $s= 1-\varepsilon \pm ix- n,\  s= -\alpha - n,\ n \in \mathbb{N}_0$.  Consequently, after straightforward calculations we express the kernel $S_\varepsilon (x,  y)$ in the form
$$S_\varepsilon (x,  y)  =\  y^\alpha \frac{\Gamma(1+\alpha) \left|\Gamma(\varepsilon-1-\alpha+ix)\right|^2}{ \Gamma(1+2\alpha)\Gamma(\varepsilon- 1/2-\alpha)} $$
$$\times \  {}_2F_3 \left( 1+\alpha,\ {3\over 2}+\alpha-\varepsilon;\ 1+ 2\alpha,\ 2+\alpha -\varepsilon- ix,\  2+\alpha -\varepsilon+  ix ;\ y \right) $$
$$+ { y^{\varepsilon-1}  \over \sqrt\pi} {\rm Re } \left[ \left({y\over 4}\right)^{ix} \ \Gamma(1- \varepsilon+\alpha -ix)  \Gamma(\varepsilon+ ix) \Gamma(-ix) \right. $$
$$\left. \times  \  {}_2F_3 \left( \varepsilon+ix,\  {1\over 2}+ ix ;\ 1+ 2 ix,\  \varepsilon- \alpha + ix,\  \varepsilon +\alpha +  ix ;\ y \right) \right] .$$
Substituting this value in (3.9) and using the reduction formula for the gamma-function, we end up with the inversion formula (3.4), completing the proof.

 \end{proof}

\section{Initial   value problem}

In this section the index transform (1.2) is employed  to investigate  the  solvability  of an initial value  problem  for the following fourth    order partial differential  equation, involving the Laplacian
$$\left[ \left(  x {\partial  \over \partial x}  + y {\partial   \over \partial y}  + 2 \right)^2+ \alpha^2 \right]   \Delta u -  {\sqrt{ x^2+y^2} + 2\alpha^2\over x^2+ y^2} \left[ x {\partial  \over \partial x}  + y {\partial   \over \partial y}\right]^2 u $$
$$ -  {3\over 2 \sqrt{ x^2+y^2} } \left[ x {\partial u  \over \partial x}  + y {\partial  u  \over \partial y}\right]  +  {u\over 2 \sqrt{ x^2+y^2}}   =0, \  (x,y) \in   \mathbb{R}^2 \backslash \{0\},\eqno(4.1)$$ 
where $\Delta = {\partial^2 \over \partial x^2} +  {\partial^2 \over \partial y ^2}$ is the Laplacian in $\mathbb{R}^2$.   In fact, writing  (4.1) in polar coordinates $(r,\theta)$, we end up with the equation  
$$ r^3 { \partial^4 u \over \partial r^4} +  r\  { \partial^4 u \over \partial r^2 \partial \theta ^2}  +  6 r^2\  { \partial^3 u \over \partial r^3}    +    { \partial^3 u \over \partial r \  \partial \theta^2} + r\ (7- \alpha^2- r)  { \partial^2 u \over \partial r^2 } $$

$$+  {\alpha^2 \over r} \  { \partial^2 u \over \partial \theta^2} +   \left( 1- \alpha^2- {5\over 2r} \right) \  { \partial u \over \partial r}  +  {u \over 2}   = 0.\eqno(4.2)$$

{\bf Lemma 2.} {\it  Let $\alpha > - 1/4,  g(\tau)  \in L_1\left(\mathbb{R}; e^{ (\beta-\pi)   |\tau|} d\tau\right),\  \beta \in (0, 2\pi)$. Then  the function
$$u_\alpha (r,\theta)=  2\sqrt\pi \int_{-\infty}^\infty  e^{\theta \tau} \ {\rm Re} \left[ K_{\alpha + i\tau}(\sqrt r)\  I_{\alpha - i\tau}(\sqrt r)\right] {g(\tau)\over \cosh(\pi\tau)} d\tau\eqno(4.3)$$
 satisfies   the partial  differential  equation $(4.2)$ on the wedge  $(r,\theta): r   >0, \  0\le \theta <  \beta$, vanishing at infinity.}

\begin{proof} The proof  is straightforward by  substitution (4.3) into (4.2) and the use of  (1.11).  The necessary  differentiation  with respect to $r$ and $\theta$ under the integral sign is allowed via the absolute and uniform convergence, which can be verified  using  inequality (3.1)  and the integrability condition $g \in L_1\left(\mathbb{R}; e^{ (\beta-\pi)  |\tau|} d\tau\right),\  \beta \in (0, 2\pi)$ of the lemma.  Finally,  the condition $ u(r,\theta) \to 0,\ r \to \infty$  is again due to (3.1). 
\end{proof}

We are ready to  formulate the initial value problem for equation (4.2) and give its solution.

{\bf Theorem 5.} {\it Let  $g(x)$ be given by formula $(3.4)$ and its transform $(G_\alpha g) (t)\equiv G_\alpha (t)$ satisfies conditions of Theorem 3.  Then  $u (r,\theta),\   r >0,  \  0\le \theta < \beta$ by formula $(4.3)$  will be a solution  of the initial value problem for the partial differential  equation $(4.2)$ subject to the initial condition}
$$u_\alpha (r,0) = G _\alpha (r).$$

Finally we will pay our attention to the so-called generalized Lebedev index transform recently considered by the author \cite{lebyak}, which contains  the square modulus of the Macdonald function as the kernel
 $$\Psi_\alpha (r,\theta) =   \int_{-\infty}^\infty e^{\theta\tau} \left| K_{\alpha+ i\tau}(\sqrt r)\right|^2  g(\tau)d\tau,\eqno(4.4)$$
where $\alpha \in \mathbb{R},\   r >0, \ 0\le \theta \le 2\pi$.  Namely, we will show  that the kernel  $ \left| K_{\alpha+ i\tau}(\sqrt r)\right|^2  $ satisfies differential equation (1.11) and, correspondingly,  the index transform (4.4) is a solution of the PDE (4.2) under the boundary condition
$$\lim_{r\to \infty} \Psi_\alpha (r,\theta) = 0.\eqno(4.5)$$

{\bf Lemma 3}. {\it Let $\alpha,\ \tau \in \mathbb{R}$.  The kernel  $\left| K_{\alpha+ i\tau}(\sqrt r)\right|^2$ is a fundamental solution of  the  differential equation $(1.11)$. }

\begin{proof}    Taking the integral representation for  the kernel in terms of the Mellin-Barnes integral via relation (8.4.23.31) in \cite{prud}, Vol. III, we have 
$$ \left| K_{\alpha+ i\tau}(\sqrt r)\right|^2 =  {1 \over 4 i \sqrt \pi } \int_{\gamma-i\infty}^{\gamma+i\infty} \Gamma(s+i\tau) \Gamma (s-i\tau) \frac{\Gamma (s+\alpha) \Gamma (s-\alpha)}
{ \Gamma(1/2+ s)  \Gamma(s)} r^{-s} ds, \eqno(4.7)$$
where $\gamma > |\alpha|.$   Hence, since the integrand in (4.7) behaves at infinity  as
$$\Gamma(s+i\tau) \Gamma (s-i\tau) \frac{\Gamma (s+\alpha) \Gamma (s-\alpha)}
{ \Gamma(1/2+s)  \Gamma(s)} = O\left( e^{-\pi |t|} |t|^{ 2\gamma - 3/2}\right), \ s= \gamma +it, \ \quad |t| \to \infty,$$
the repeated differentiation with respect to $r$ under the integral sign  is permitted.  Therefore following   the same scheme as in the proof of Lemma 1, we derive  
$$ \left( r{d\over dr} \right)^2 \   \left| K_{\alpha+ i\tau}(\sqrt r)\right|^2 =  {1\over 4i \sqrt \pi } \int_{\gamma-i\infty}^{\gamma+i\infty} \Gamma(s+i\tau) \Gamma (s-i\tau) \frac{ s^2  \  \Gamma (s+\alpha) \Gamma (s-\alpha)}  { \Gamma(1/2+s)  \Gamma(s)} r^{-s} ds$$
$$= -  \tau^2 \left| K_{\alpha+ i\tau}(\sqrt r)\right|^2    + {1\over 4i \sqrt \pi } \int_{\gamma-i\infty}^{\gamma+i\infty} \Gamma(1+ s+i\tau) \Gamma (1+ s-i\tau) \frac{ \Gamma (s+\alpha) \Gamma (s-\alpha)}  { \Gamma(1/2+s)  \Gamma(s)} r^{-s} ds$$
  $$= - \tau^2 \left| K_{\alpha+ i\tau}(\sqrt r)\right|^2    + {1\over 4 i\sqrt \pi } \int_{1+\gamma-i\infty}^{1+\gamma+i\infty} \Gamma( s+i\tau) \Gamma ( s-i\tau) \frac{ \Gamma (s-1 +\alpha) \Gamma (s-1 -\alpha)}  { \Gamma(s- 1/2)  \Gamma(s-1)} r^{1-s} ds$$
  $$= -  \tau^2 \left| K_{\alpha+ i\tau}(\sqrt r)\right|^2    + {1\over 4i \sqrt \pi } \int_{1+\gamma-i\infty}^{1+\gamma+i\infty} \Gamma( s+i\tau) \Gamma ( s-i\tau) \frac{ (s- 1/2) (s-1) \Gamma (s +\alpha) \Gamma (s -\alpha)}  {(s-1 +\alpha) (s-1 -\alpha)  \Gamma(s+ 1/2)  \Gamma(s)} r^{1-s} ds.$$
  Hence, shifting  the contour to the left by Cauchy's theorem, we multiply the latter equality by $r^{-\alpha} $ and differentiate it with respect to $r$ again. Then  
$${d\over dr} \left[ r^{-\alpha}  \left( r{d\over dr} \right)^2 \   \left| K_{\alpha+ i\tau}(\sqrt r)\right|^2 \right] =
-  \tau^2  {d\over dr} \left[ r^{-\alpha}  \left| K_{\alpha+ i\tau}(\sqrt r)\right|^2 \right] $$
$$-  {1\over 4i \sqrt \pi} \int_{1+\gamma-i\infty}^{1+\gamma+i\infty} \Gamma( s+i\tau) \Gamma ( s-i\tau) \frac{ (s- 1/2) (s-1) \Gamma (s +\alpha) \Gamma (s -\alpha)}  { (s-1 -\alpha)  \Gamma(s+ 1/2)  \Gamma(s)} r^{-s-\alpha } ds.$$
In the same fashion we deduce the equality 
$$ {d\over dr} \left[ r^{1+2 \alpha}  {d\over dr} \left[ r^{-\alpha}  \left( r{d\over dr} \right)^2 \   \left| K_{\alpha+ i\tau}(\sqrt r)\right|^2 \right]  \right] = -  \tau^2 {d\over dr} \left[ r^{1+2 \alpha} {d\over dr} \left[ r^{-\alpha}  \left| K_{\alpha+ i\tau}(\sqrt r)\right|^2 \right] \right] $$
$$+  {1\over 4i \sqrt \pi} \int_{1+\gamma-i\infty}^{1+\gamma+i\infty} \Gamma( s+i\tau) \Gamma ( s-i\tau) \frac{ (s- 1/2) (s-1) \Gamma (s +\alpha) \Gamma (s -\alpha)}  { \Gamma(s+ 1/2)  \Gamma(s)} r^{-s + \alpha } ds.$$
Consequently,  minding (4.7), we have 
$$ r^{-\alpha}  {d\over dr} \left[ r^{1+2 \alpha}  {d\over dr} \left[ r^{-\alpha}  \left( r{d\over dr} \right)^2 \   \left| K_{\alpha+ i\tau}(\sqrt r)\right|^2 \right]  \right] = -  \tau^2  r^{-\alpha} {d\over dr} \left[ r^{1+2 \alpha} {d\over dr} \left[ r^{-\alpha}  \left| K_{\alpha+ i\tau}(\sqrt r)\right|^2 \right] \right] $$
$$-  {1\over 4i \sqrt \pi} {d\over dr} \int_{1+\gamma-i\infty}^{1+\gamma+i\infty} \Gamma( s+i\tau) \Gamma ( s-i\tau) \frac{ (s- 1/2)  \Gamma (s +\alpha) \Gamma (s -\alpha)}  { \Gamma(s+ 1/2)  \Gamma(s)} r^{1-s} ds$$
$$= -  \tau^2  r^{-\alpha} {d\over dr} \left[ r^{1+2 \alpha} {d\over dr} \left[ r^{-\alpha}  \left| K_{\alpha+ i\tau}(\sqrt r)\right|^2 \right] \right]  + {1\over 2} {d\over dr} \left[ r   \left| K_{\alpha+ i\tau}(\sqrt r)\right|^2 \right]  +   {d\over dr}   r^2    {d\over dr} \left| K_{\alpha+ i\tau}(\sqrt r)\right|^2.$$
Thus fulfilling the differentiation, we come again to (1.11). 

\end{proof}

{\bf Theorem 6.} {\it  Let $\alpha  \in \mathbb{R},  g(\tau)  \in L_1\left(\mathbb{R}; e^{ (\beta-\pi)   |\tau|} d\tau\right),\  \beta \in (0, 2\pi)$. Then  the function $(4.4)$  satisfies   the partial  differential  equation $(4.2)$ on the wedge  $(r,\theta): r   >0, \  0\le \theta <  \beta$, vanishing at infinity.}

\bigskip
\centerline{{\bf Acknowledgments}}
\bigskip

The work was partially supported by CMUP (UID/MAT/00144/2013), which is funded by FCT(Portugal) with national (MEC) and European structural funds through the programs FEDER, under the partnership agreement PT2020.

\bibliographystyle{amsplain}

\end{document}